\nonstopmode
\documentclass[]{amsart}
\usepackage{amssymb}
\usepackage{graphicx}
\usepackage{epsfig}
\usepackage[all]{xy}
\newdir{ >}{{}*!/-10pt/@{>}}
\allowdisplaybreaks

\def\undersetbrace#1\to#2{\underbrace{#2}_{#1}}
\def\oversetbrace#1\to#2{\overbrace{#2}^{#1}}
\def\AMSunderset#1\to#2{\underset{#1}{#2}}
\def\AMSoverset#1\to#2{\overset{#1}{#2}}

\def\East#1#2{\overset{#1}{\longrightarrow}}
\swapnumbers

\newtheorem{prop}[subsection]{Proposition}
\newtheorem*{prop*}{Proposition}

\newtheorem*{thm*}{Theorem}
\newtheorem{lem}[subsection]{Lemma}
\newtheorem*{lem*}{Lemma}
\newtheorem{cor}[subsection]{Corollary}
\newtheorem*{cor*}{Corollary}

\newtheorem*{conj*}{Conjecture}
\newenvironment{demo}[1]{\par\smallskip\noindent{\bf #1.}}{\par\smallskip}

\def\cit#1#2{\ifx#1!\cite{#2}\else#2\fi} 
\def\ign#1{}             
\def\o{\circ}

\def\al{\alpha}
\def\be{\beta}

\def\th{\theta}

\def\ph{\varphi}

\def\Ph{\Phi}

\def\i{^{-1}}
\def\x{\times}

\let\on=\operatorname

\def\AMSonly#1{}
\begin{document}
\title[]
{The homotopy type of the space of degree 0 immersed plane curves
}
\author{Hiroki Kodama, Peter W. Michor}
\address{
Hiroki Kodama:
Erwin Schr\"odinger Institut f\"ur Mathematische Physik,
Boltzmanngasse 9, A-1090 Wien, Austria
}
\email{kodama@ms.u-tokyo.ac.jp}
\address{
Peter W. Michor:
Fakult\"at f\"ur Mathematik, Universit\"at Wien,
Nordbergstrasse 15, A-1090 Wien, Austria; {\it and:}
Erwin Schr\"odinger Institut f\"ur Mathematische Physik,
Boltzmanngasse 9, A-1090 Wien, Austria
}
\email{Peter.Michor@univie.ac.at}
\date{\today}
\thanks{PWM was supported by FWF Project P~17108. 
HK was supported by a JSPS-IH\'ES-EPDI fellowship}

\begin{abstract} The space 
$B^0_i=\on{Imm}^0(S^1,\mathbb R^2)/\on{Diff}(S^1)$ of
all immersions of rotation degree 0 in the plane modulo reparametrizations
has homotopy groups $\pi_1(B^0_i)=\mathbb Z$, $\pi_2(B^0_i)=\mathbb Z$, and 
$\pi_k(B^0_i)=0$ for $k\ge 3$.
\end{abstract}
\def\LaTeXonly{}

\maketitle

\section{Introduction} \label{nmb:1}
For an immersion $c:S^1\to \mathbb R^2$ the (rotation) degree is the winding number of
$c':S^1 \to \mathbb R^2$ around 0. Let $\on{Imm}^k(S^1,\mathbb R^2)$ denote the
connected smooth Fr\'echet manifold of all immersions of degree $k$.
It was shown in \cite{MM1} that for $k\ne 0$ the space $\on{Imm}^k(S^1,\mathbb R^2)$
contains a copy of $S^1$ as a smooth strong deformation retract and that
the infinite dimensional orbifold $\on{Imm}^k(S^1,\mathbb R^2)/\on{Diff}^+(S^1)$
is contractible, where $\on{Diff}^+(S^1)$ denotes the regular Fr\'echet Lie
group of orientation preserving diffeomorphisms. The proof in \cite{MM1} 
consists in expanding the classical proof of the theorem of Whithney and Graustein 
(see \cite{G1}, \cite{G2}, \cite{Wh})
into the construction of an $S^1$-equivariant smooth deformation
retraction. For $k=0$ this did not work. 

In this paper we treat the case $k=0$. 
In section \ref{nmb:2} we first give simple argument which shows that 
$\pi_1(B^0(S^1,\mathbb R^2))$
contains $\mathbb Z$. In section \ref{nmb:3} we give a more involved proof that 
$\on{Imm}^0(S^1,\mathbb R^2)$ is homotopy equivalent to $S^1$. 
In section \ref{nmb:4} we show that factoring out $\on{Diff}^+(S^1)$ gives
a fibration with homotopically trivially embedded fiber, 
and then the homotopy sequence shows that 
$\pi_1(B^{0,+}(S^1,\mathbb R^2))=\mathbb Z$, 
$\pi_2(B^{0,+}(S^1,\mathbb R^2))=\mathbb Z$, and
$\pi_k(B^{0,+}(S^1,\mathbb R^2))=0$ for $k>2$.
Factoring out the larger group $\on{Diff}(S^1)$ gives a two-sheeted
covering and the final result.

PM thanks the audience of the lecture course SS~2005 and Peter Zvengrowksi
for helpful discussions.
HK thanks the Erwin Schr\"odinger Institute for the warm
hospitality during his stay in Vienna.

\section{A simple proof that  $\mathbb Z\subseteq\pi_1(B^0(S^1,\mathbb R^2))$} 
\label{nmb:2}

\begin{prop}
$\on{Imm}^0(S^1,\mathbb R^2) /\on{Diff}^+(S^1)$ is not contractible.
\end{prop}

\begin{demo}{Proof}
We shall view a curve $c\in\on{Imm}(S^1,\mathbb R^2)$
as a $2\pi$-periodic plane valued function.
A smooth function $ a = a(c, \quad) : \mathbb R \to \mathbb R$
is called an argument of a curve $c$ if
$$
\frac{c'(\th)}{\lvert c'(\th) \rvert} = \exp (i\,a(\th));
$$
it is unique up to addition of an inter multiple of $2\pi$.
If the curve $c$ has degree $k$ then $a(\th+2\pi)-a(\th) = 2k\pi$.
Thus, a curve $c$ is in $\on{Imm}^0(S^1,\mathbb R^2)$
if and only if some (any) argument of $c$ is $2\pi$-periodic.
For a curve $c\in\on{Imm}^0(S^1,\mathbb R^2)$, 
we define the \textit{average argument}
$\al(c) \in S^1$ by
$$
\al(c) = \exp 
\left(
\frac{i}{l(c)} \int_{0}^{2\pi} a(c,\th) \lvert c'(\th) \rvert d\th \right),
$$
which does not depend on the choice of $a(c,\quad)$ and defines a 
well-defined smooth mapping
$\al:\on{Imm}^0(S^1,\mathbb R^2)\to S^1$.
Also, since any argument $a$ of a degree 0 curve is $2\pi$-periodic,
$\al(c)$ is invariant under the action of
$\on{Diff}^+(S^1)$.   
So we can view $\al$ as a map 
$$
\al\colon B^{0,+}(S^1,\mathbb R^2)= \on{Imm}^0(S^1,\mathbb R^2) /\on{Diff}^+(S^1) \to
S^1.
$$
For $\ph \in S^1\subset \mathbb C=\mathbb R^2$, the rotation map 
$\ph:\mathbb R^2\to\mathbb R^2$ act on 
$B^{0,+}(S^1,\mathbb R^2)$ and 
obviously
$$
\al(\ph.c) = \ph.\al(c) .
$$
So choosing a free orbit $S^1.C$ for the rotation action of $S^1$ on 
$B^{0,+}(S^1,\mathbb R^2)$, the composition
$$
S^1.C\hookrightarrow  \on{Imm}^0(S^1,\mathbb R^2) /\on{Diff}^+(S^1)
\East{\al}{} S^1
$$
equals the identity on $S^1$, thus 
$\pi_1(S^1)=\mathbb Z\subset \pi_1(B^{0,+}(S^1,\mathbb R^2))$.

Moreover, $\al(c(-\quad))=-\al(c)$ implies that $\al$ factors as follows,
where the vertical arrows are 2-sheeted coverings:
$$
\xymatrix{
B^{0,+}(S^1,\mathbb R^2) \ar@{=}[r] 
&\on{Imm}^0(S^1,\mathbb R^2) /\on{Diff}^+(S^1) \ar[r]^{\qquad\qquad\al} \ar[d]_2 
& S^1 \ar[d]^2 \\
B^0(S^1,\mathbb R^2) \ar@{=}[r]
&\on{Imm}^0(S^1,\mathbb R^2) /\on{Diff}(S^1) \ar[r]^{\qquad\qquad\bar\al} 
& S^1 \\
}
$$
Thus we also get in a similar way
$\pi_0(S^1)=\mathbb Z\subset \pi_0(B^0(S^1,\mathbb R^2))$.
\qed\end{demo}

\section{The homotopy type of $\on{Imm}^0(S^1,\mathbb R^2)$}
\label{nmb:3}

\begin{prop}\label{nmb:3.1}
The space $\on{Imm}^{0}(S^1,\mathbb R^2)$ of degree 0 immersions in the
plane is homotopy equivalent to $S^1$. 
\end{prop}

\begin{demo}{Proof}
This will follow from \ref{nmb:3.2} -- \ref{nmb:3.5} below.
\qed\end{demo}

\subsection{}\label{nmb:3.2}
Let  
$\on{Imm}^{0,*}(S^1,\mathbb R^2) := \{c\in\on{Imm}^0(S^1,\mathbb R^2) ;\;c(0)=0\}$. 
Clearly we have 
$\on{Imm}^0(S^1,\mathbb R^2) 
\cong \on{Imm}^{0,*}(S^1,\mathbb R^2) \times \mathbb R^2$
and
$\on{Imm}^0(S^1,\mathbb R^2) \sim \on{Imm}^{0,*}(S^1,\mathbb R^2)$,
where $\cong$ denotes homoeomorphism and $\sim$ homotopy equivalent.
Let us define a map
\begin{align*}
&\Ph \colon \on{Imm}^{0,*} \to 
C^\infty(S^1,\mathbb R_+) \times C^\infty(S^1,S^1)
\\&
\Ph(c)(\theta) = 
\left(\lvert c_\theta(\theta)\rvert,
\frac{c_\theta(\theta)}{\lvert c_\theta(\theta)\rvert}
\right)  =: (v(\theta),e(\theta)).
\end{align*}
The map $\Ph$ is injective.
For $(v,e) = \Ph(c)$, the winding number of $e$ equals the degree $0$ of
$c$ and thus
$\int_0^{2\pi}v.e\,d\theta=0$.

\begin{lem*}
The length of the image of $e$ is greater than $\pi$.
\end{lem*}

\begin{demo}{Proof}
If not, there exists a number $r \in \mathbb R$ such that
$$
\exp(ir) \in \on{Im}(e) 
\subset \exp(i[r-\pi/2,r+\pi/2]).
$$
Then, $\langle \exp(ir),e(\theta) \rangle$ is nonnegative for any $\theta$
and strictly positive for some $\theta$.
Therefore
$\int_0^{2\pi} \langle \exp(ir), v.e \rangle \,d\theta>0$.
This contradicts 
$$
\int_0^{2\pi} \langle \exp(ir), v.e \rangle \,d\theta
= \Big\langle \exp(ir), \int_0^{2\pi}  v.e \,d\theta \Big\rangle
=\langle \exp(ir), 0 \rangle
= 0. \qed
$$
\end{demo}

\subsection{}\label{nmb:3.3}
Let us define the set
$$
C^{\infty,0}_{>\pi}(S^1,S^1)
= \{ e \in C^\infty(S^1,S^1) ;\; \on{deg}(e) =0,\, 
\on{length}({\on{Im}(e)}) > \pi \}
$$
and consider the map
$$
\on{pr}_2 \circ \Ph \colon
\on{Imm}^{0,*}(S^1,S^1) \to C^{\infty,0}_{>\pi}(S^1,S^1),
$$
where $\on{pr}_2$ denotes the second projection.

\begin{lem*}
The map 
$
\on{pr}_2 \circ \Ph \colon
\on{Imm}^{0,*}(S^1,S^1) \to C^{\infty,0}_{>\pi}(S^1,S^1),
$
is surjective, has contractible fibers, admits a global smooth section, and
is a homotopy equivalence.
\end{lem*}

\begin{demo}{Proof}
For a map $e \in C^{\infty,0}_{>\pi}(S^1,S^1)$,
there exist points $\theta_1$, $\theta_2$, $\theta_3$
such that $0 \in \on{int}([e(\theta_1),e(\theta_2),e(\theta_3)])$,
where $[\,\cdot,\cdot,\cdot\,]$ denotes the convex hull of three points.
Let $v_1\in C^\infty(S^1,\mathbb R_{>0})$ be a map such that
$\int_0^{2\pi} v_1 d\theta = 1$ and
$v_1(\theta)$ is close to $0$ if 
$\theta$ is not close to $\theta_1$.
Then $\int_0^{2\pi}v_1.e\,d\theta$ is close to $e(\theta_1)$.
We also define $v_2$ and $v_3$ similarly, so that
$$
0 \in \on{int}\left(\left[\int_0^{2\pi}v_1.e\,d\theta,
\int_0^{2\pi}v_2.e\,d\theta, \int_0^{2\pi}v_3.e\,d\theta \right] \right).
$$
Therefore there exist positive numbers $a_1$, $a_2$, $a_3$
with  
$$
a_1\int_0^{2\pi}v_1.e\,d\theta+ a_2\int_0^{2\pi}v_2.e\,d\theta+
a_3\int_0^{2\pi}v_3.e\,d\theta=0.
$$
Define $c$ by
$$
c(\theta) = 
\int_0^\theta  (a_1v_1(u)+a_2v_2(u)+a_3v_3(u)) e(u) \,du.
$$
Then $c$ is in $\on{Imm}^{0,*}$ and 
$(\on{pr}_2 \circ \Ph) (c) = e$, which means that 
$\on{pr}_2 \circ \Ph$ is surjective.

We next show that for any $e \in C^{\infty,0}_{>\pi}(S^1,S^1)$,
the inverse image $(\on{pr}_2 \circ \Ph)^{-1}(e)$ is contractible.
Namely, let
$V(e)\subset C^\infty(S^1,\mathbb R_+)$ be given by
$$
V(e) = 
\left\{ 
v \in C^\infty(S^1,\mathbb R_+) ;\; \int_0^{2\pi} v.e \,d\theta =0
\right\},
$$
an open convex subset of the linear subspace
$\{v \in C^\infty(S^1,\mathbb R_+) ;\; \int_0^{2\pi} v.e \,d\theta =0\}\subset
C^\infty(S^1,\mathbb R)$. 
Thus $V(e)$ is contractible for each $e$.
Moreover, $V(e)$ is homeomorphic to $(\on{pr}_2 \circ \Ph)\i(e)$
by the map $\on{pr}_1 \circ \Ph \colon (p_2 \circ \Ph)^{-1}(e) \to V(e)$.

For
fixed $\theta_1,\theta_2,\theta_3$ the construction above works for each
$e\in  C^{\infty,0}_{>\pi}(S^1,S^1)$ for which $0$ is contained in the
interior of the convex hull of $e(\theta_1),e(\theta_2), e(\theta_3)$;
these $e$ form an open set in $C^{\infty,0}_{>\pi}(S^1,S^1)$ on which we
get a continuous (even smooth) section of $p_2\circ \Ph$. Open sets like that
cover $C^{\infty,0}_{>\pi}(S^1,S^1)$. So we get smooth local sections whose
domains cover the base. Since the base is open in a nuclear Fr\'echet
space, it is smoothly paracompact (see \cite{KM},~16.10) we can use convexity
of all fibers and a smooth partition of unity on the base 
$C^{\infty,0}_{>\pi}(S^1,S^1)$ to construct a global smooth section $s$. 

Finally, since all fibers are convex, there is a smooth strong fiber
preserving deformation retraction of $\on{Imm}^{0,*}(S^1,S^1)$ onto the
image the global section $s$.
\qed
\end{demo}

\subsection{}\label{nmb:3.4}
To study the topology of $C^{\infty,0}_{>\pi}(S^1,S^1)$,
we introduce the set of $2\pi$-periodic functions
$$
C^{\infty,p}(\mathbb R,\mathbb R) 
= \{ c \in C^{\infty}(\mathbb R,\mathbb R) ;\; c(\theta+2\pi) = c(\theta) \}
.$$
For $c \in C^{\infty,p}(\mathbb R,\mathbb R)$ let $\on{Var}(c) = \max c - \min c$ and
let $\on{Ave}(c) = \frac{1}{2\pi} \int_0^{2\pi} c\,d\theta$.
For $k \geq 0$, 
$
C^{\infty,p}_{>k}(\mathbb R,\mathbb R) 
= \{ c \in C^{\infty,p}(\mathbb R,\mathbb R) ;\; \on{Var}(c)>k \}
$.
Define a diffeomorphism $g \colon C^{\infty,p}_{>0}(\mathbb R,\mathbb R) 
\to C^{\infty,p}_{>\pi}(\mathbb R,\mathbb R)$ by
$$
g(c) = 
\frac{\on{Var}(c)+\pi}{\on{Var}(c)}
(c-\on{Ave}(c))+\on{Ave}(c).
$$
The diffeomorphism $g$ satisfies $g(c(\quad+2n\pi))=g(c)(\quad+2n\pi)$,
thus induces the diffeomorphism
$$
\tilde{g} \colon 
C^{\infty,0}_{>0}(S^1,S^1) \to
C^{\infty,0}_{>\pi}(S^1,S^1) 
,$$
where $C^{\infty,0}_{>0}(S^1,S^1)$ denotes the set of
nonconstant smooth maps of degree 0 in $C^{\infty}(S^1,S^1)$.

\subsection{}\label{nmb:3.5}
We consider now the evaluation $\on{ev}_1$ at $1\in S^1$ whose fiber at
$1\in S^1$ is the smooth manifold of based smooth loops of degree 0 in
$S^1$, with the constant loop $1$ deleted:
\begin{equation*}
C^{\infty,0}((S^1,1),(S^1,1))\setminus \{1\} \hookrightarrow
C^{\infty,0}_{>0}(S^1,S^1) \East{\on{ev}_1}{} S^1
\tag{1}\end{equation*}

\begin{lem*}
The map
$\on{ev}_1:C^{\infty,0}_{>0}(S^1,S^1)\to S^1$ is a smooth trivial
fibration with a global section and smoothly contractible fibers. Moreover,
it is a homotopy equivalence.
\end{lem*}

\begin{demo}{Proof}
A smooth section $s:S^1\to C^{\infty,0}_{>0}(S^1,S^1)$ of $\on{ev}_1$ is given by 
$s(\ph)(\th)=\ph.\exp(i\on{Im}(\th))$. The fiber of $\on{ev}_1$ over $\ph$
is the space $C^{\infty,0}((S^1,1),(S^1,\ph))\setminus \{\ph\}$ consisting
of all non-constant smooth loops of degree 0 mapping $1$ to $\ph$, which is
diffeomorphic to the fiber $C^{\infty,0}((S^1,1),(S^1,1))\setminus \{1\}$
via multiplication by $\ph$. 

It remains to show that the fiber $C^{\infty,0}((S^1,1),(S^1,1))\setminus \{1\}$
is contractible. Via lifting to the univeral cover, 
$C^{\infty,0}((S^1,1),(S^1,1))$ is diffeomorphic to the space 
$\{f\in C^{\infty,p}(\mathbb R,\mathbb R): f(0)=0\}$ of periodic functions
mapping 0 to 0. Via Fourier expansion $f(t)=\sum_{n\in\mathbb Z}a_n\exp(int)$ 
this is isomorpic to the space of all
rapidly decreasing complex sequences $(a_k)_{k\in\mathbb Z}$ with 
$\overline{a_k}=a_{-k}$ and $\sum_k a_k=0$. This space is
isomorphic to the space $\mathfrak s$ of rapidly decreasing sequences $(b_n)_{n\ge 1}$
by $a_n=b_n$ for $n\ge 1$, $a_{-n}=\overline b_n$, and
$a_0=2\on{Re}(\sum_{n\ge1} b_n)$.

Now we have to show that this is still contractible if we
remove the constant sequence 0. Then it is homotopy equivalent to its
intersection with the 
sphere in $\ell^2$, i.e., to the space 
$S:=\{b\in\mathfrak s: \sum_{n\ge1} b_n^2=1\}$.
But this is contractible by a standard argument which is explained on page
513 of \cite{KM} for the space of finite sequences. 
Namely, 
consider the homotopy 
$A:\mathfrak s \x [0,1]\to\mathfrak  s $ through isometries which 
is given by $A_0=\on{Id}$ and by 
\begin{multline*}
A_t(b_1,b_2,b_3,\ldots)=(b_1,\dots,b_{n-2}, b_{n-1}\cos\th_n(t), 
     b_{n-1}\sin\th_n(t),\\
b_n\cos\th_n(t),b_n\sin\th_n(t),
     b_{n+1}\cos\th_n(t),b_{n+1}\sin\th_n(t),\ldots)
\end{multline*}
for $\frac1{n+1}\leq t\leq\frac1n$, where 
$\th_n(t)=\ph(n((n+1)t-1))\frac\pi2$ for a fixed smooth function 
$\ph:\mathbb R\to\mathbb R$ which is 0 on $(-\infty,0]$, grows 
monotonely to 1 in $[0,1]$, and equals 1 on $[1,\infty)$.
The mapping $A$ is Lipschitz continuous for each seminorm
$\|b\|_k=\sup\{|b_n|n^k:n\ge 1\}$ of $\mathfrak s$
with constant $2^k$, and is isometric for
$\ell^2$. 
Then       
$A_{1/2}(b_1,b_2,\ldots)=(b_1,0,b_2,0,\ldots)$ is in 
$\mathfrak s _{\text{odd}}$, and on the other hand 
$A_1(b_1,b_2,\ldots)=(0,b_1,0,b_2,0,\ldots)$ is in 
$\mathfrak s _{\text{even}}$.
This is a variant of a homotopy constructed by \cite{Ra}.
Now $A_t|S$ for $0\le t\le 1/2$ is a homotopy on 
$S$ between the identity and 
$A_{1/2}(S)\subset \mathfrak s _{\text{odd}}$. The latter 
set is contractible, for example in a stereographic chart. 
\qed\end{demo}

\subsection{}\label{nmb:3.6}
If we put together all mappings constructed above we get the following
commutative diagram where we indicate isomorphism $\cong$, homotopy
equivalence $\sim$, or 2-sheeted covering 2, and a free orbit $S^1.c$ for
the rotation action on $\on{Imm}^0$:
$$
\xymatrix{
& S^1 \ar[r]^{=} & S^1 \ar[r]^2 & S^1 \\
S^1.c \ar[r]^{\subset} \ar[ur]^{\cong} 
& \on{Imm^0} \ar[u]_{\al} \ar@{->>}[r] \ar[d]^{\tilde g\i\o\on{pr}_2\o\Ph}_{\sim}
& B^{0,+} \ar[u]^{\al} \ar[r]^2
& B^0 \ar[u]^{\bar\al} \\
S^1 \ar[u]^{=} 
& C^{\infty,0}_{>0} \ar[l]^{\sim}_{\on{ev}_1} & &
}
$$

\section{The homotopy type of $B^0(S^1,\mathbb R^2)$ }\label{nmb:4}

\begin{prop}\label{nmb:4.1}
The mapping $\on{Imm}^0(S^1,\mathbb R^2)\to B^{0,+}(S^1,\mathbb R^2)$ is a
(Serre) fibration. 
\end{prop}

\begin{demo}{Proof}
First we replace $\on{Imm}^0(S^1,\mathbb R^2)$ by 
the subset 
$\on{Imm}^0_a(S^1,\mathbb R^2)$ consisting of all immersions which
are parametrized by scaled arc-length which is
a strong deformation retract, see \cite{MM1},~2.6. The normalizer
of the $\on{Diff}^+(S^1)$-action on it is just the action of $S^1$ which
shifts the initial point. 
We have to show that for any compactly generated space $P$ and a homotopy 
$h:[0,1]\x P\to B^{0,+}$ whose initial value $h(0,\quad)$ admits a
continuous lift there exists a continuous lift of the whole homotopy:
$$
\xymatrix{
\{0\}\x P \ar[rr]^{H(0,\quad)} \ar[d]^{\subset} && 
                                     \on{Imm}^0_a(S^1,\mathbb R^2) \ar[d] \\
[0,1]\x P \ar[rr]^{h} \ar@{-->}[urr]^H &&       B^{0,+} (S^1,\mathbb R^2)
}
$$
To get the lift $H$ we just have to specify the initial point coherently from
$H(0,p)(1)$ over $[0,1]\ni t\mapsto h(t,p)$. 

For that we need a description of the elements in  
$B^{0,+} (S^1,\mathbb R^2)$. A point $C$ in it can be described by the following
data: 

For some $n$ and $i=1,\dots,n$, there are open sets $U_i=U_i(C)\subseteq \mathbb R^2$, 
smooth functions
$f_i=f_i(C):U_i\to \mathbb R$ such that $f_i\i(0)=:C_i$ is a component $C_i$ of $C$ with
$\on{grad}(f_i)$ is a unit vector field with flow lines unit speed straight
lines passing orthogonally through $C_i$ in such a way that for $x\in C_i$
the frame consisting of $\on{grad}(f_i)(x)$ and the unit tangent to $C_i$
at $x$ is positively oriented. The unparameterized smooth oriented
1-manifolds $C_1, C_2,\dots,C_n$ (in that order) describe $C$. Note that there is a choice
for the $U_i$ and their cyclic order, but then the $f_i$ are unique.

For every $p\in P$ the initial point $H(0,p)(1)$ lies in some component
$h(0,p)_i$ of $h(0,p)$, and we may move it orthogonally along
$\on{grad}(f_i(h(t,p)))$ to get a coherent choice of initial points. This
takes care of the lift $H$.
\qed\end{demo}

\begin{lem}\label{nmb:4.2}
The fiber $\on{Diff}^+(S^1)$ maps homotopically trivial
into the fibration $\on{Imm}^0(S^1,\mathbb R^2)\to B^{0,+}(S^1,\mathbb R^2)$.
\end{lem}

\begin{demo}{Proof}
As in the proof of \ref{nmb:4.1} we consider the space
$\on{Imm}^0_a(S^1,\mathbb R^2)$ of degree 0 immersions with constant speed
parametrizations. Let $c$ be the unit speed parameterized horizontal figure
eight, and consider the diagram where $c^*(f)=c\o f$:
$$
\xymatrix{
S^1 \ar[rr]^{c^*} \ar[d]^{\subset}_{\sim} && 
    \on{Imm}^0_a  \ar@{->>}[rrd] \ar[d]^{\subset}_{\sim} && \\
\on{Diff}^+(S^1) \ar[rr]^{c^*} &&
  \on{Imm^0} \ar@{->>}[rr] \ar[d]^{\tilde g\i\o\on{pr}_2\o\Ph}_{\sim} &&
  B^{0,+} \\
S^1 && C^{\infty,0}_{>0} \ar[ll]^{\sim}_{\on{ev}_1} &&
}
$$
We have to show that the mapping from the upper left $S^1$ to the lower left 
$S^1$ is nullhomotopic. It is essentially (suppressing $\tilde g\i$) given by 
$\be\mapsto \frac{c'(\be)}{|c'(\be)|}$. From the figure
$$
\epsfig{width=8cm,file=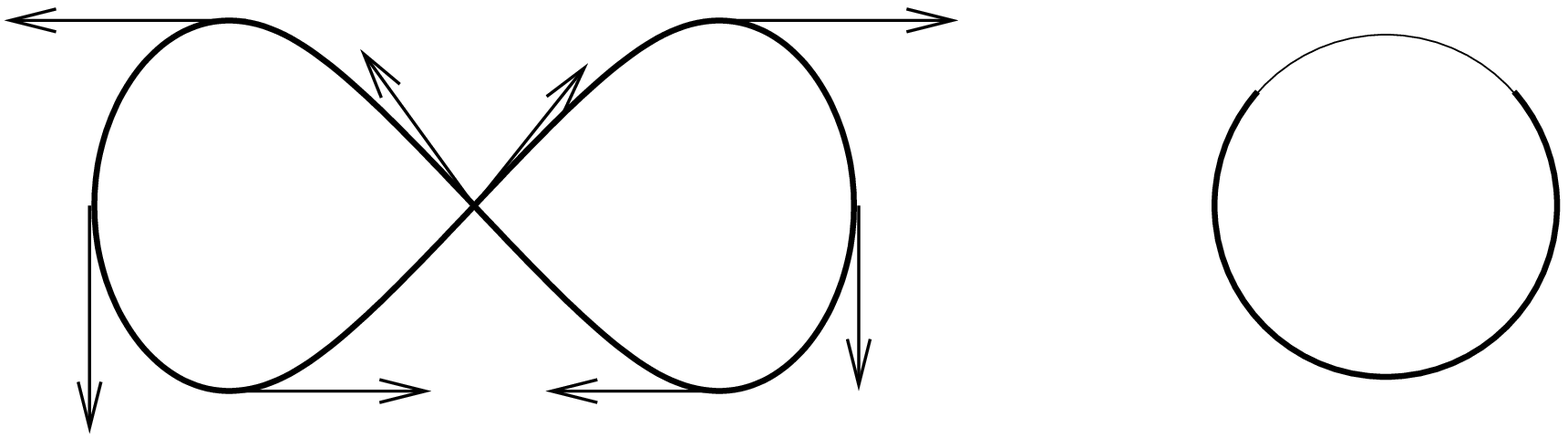}
$$
we see that this mapping covers everything below the northern polar region
twice and avoids the northern polar region, so it is nullhomotopic. 
\qed\end{demo}

\begin{cor}\label{nmb:4.3}
We have the following homotopy groups: 
\begin{align*}
\pi_1(B^{0,+}(S^1,\mathbb R^2))&=\mathbb Z, &\qquad 
  \pi_1(B^{0}(S^1,\mathbb R^2))&=\mathbb Z,
\\
\pi_2(B^{0,+}(S^1,\mathbb R^2))&=\mathbb Z, &\qquad 
  \pi_2(B^{0}(S^1,\mathbb R^2))&=\mathbb Z,
\\
\pi_k(B^{0,+}(S^1,\mathbb R^2))&=0,&\qquad	
  \pi_k(B^{0}(S^1,\mathbb R^2))&=0\quad\text{  for } k>2.
\end{align*}
\end{cor}

\begin{demo}{Proof}
By \ref{nmb:4.1} we have the long exact homotopy sequence
$$
\dots \to \pi_k(S^1) \East{0}{} \pi_k(\on{Imm}^0_a) \to \pi_k(B^{0,+}) \to
\pi_{k-1}(S^1) \to \dots
$$
and by section \ref{nmb:3} the space $\on{Imm}^0_a$ is homotopy equivalent
to $S^1$. This gives the homotopy groups of $B^{0,+}(S^1,\mathbb R^2))$.
Since $B^{0,+}(S^1,\mathbb R^2))\to B^{0}(S^1,\mathbb R^2))$ is a
two-sheeted covering, we can also read of the homotopy groups of
$B^{0,+}(S^1,\mathbb R^2))$.
\qed\end{demo}

\bibliographystyle{plain}

\end{document}